\documentclass[11pt]{article}

\usepackage{macros}
\usepackage{mathrsfs}
\usepackage[left=1.05in, right=1.05in, top=1.00in, bottom=1.05in]{geometry}

\title{Density of monochromatic infinite subgraphs}
\author{Louis DeBiasio$^{1}$, Paul McKenney$^{2}$}
\date{\today}

\begin{document}

\maketitle
\noindent\footnotetext[1]{Department of Mathematics, Miami University {\tt debiasld@miamioh.edu}. Research supported in part by Simons Foundation Collaboration Grant \# 283194.}
\noindent\footnotetext[2]{Department of Mathematics, Miami University {\tt mckennp2@miamioh.edu}.}

\begin{abstract}
For any countably infinite graph $G$, Ramsey's theorem guarantees an infinite monochromatic copy of $G$ in any $r$-coloring of the edges of the countably infinite complete graph $K_\NN$.  Taking this a step further, it is natural to wonder how ``large'' of a monochromatic copy of $G$ we can find with respect to some measure -- for instance, the density (or upper density) of the vertex set of $G$ in the positive integers.  Unlike finite Ramsey theory, where this question has been studied extensively, the analogous problem for infinite graphs has been mostly overlooked.

In one of the few results in the area, Erd\H{o}s and Galvin proved that in every 2-coloring of $K_\NN$, there exists a monochromatic path whose vertex set has upper density at least $2/3$, but it is not possible to do better than $8/9$. They also showed that for some sequence $\ep_n\to 0$, there exists a monochromatic path $P$ such that for infinitely many $n$, the set $\{1,2,...,n\}$ contains the first $(\frac{1}{3+\sqrt{3}}-\ep_n)n$ vertices of $P$, but it is not possible to do better than $2n/3$.  We improve both results, in the former case achieving an upper density at least $3/4$ and in the latter case obtaining a tight bound of $2/3$.  We also consider related problems for directed paths, trees (connected subgraphs), and a more general result which includes locally finite graphs for instance.

\end{abstract}

\section{Introduction}

Throughout the paper we let $K_\NN$ be the complete graph on the positive integers and $\vec{K}_\NN$ be the complete symmetric digraph on the positive integers.  By $r$-coloring of $G$, we mean a partition of the edges of $G$ into $r$ parts (i.e.\ an $r$-edge-coloring).  We say that a collection of subgraphs $\{H_1, H_2, \dots, H_k\}$ covers (partitions) $G$ if the collection of vertex sets $\{V(H_1), V(H_2), \dots, V(H_k)\}$ forms a cover (partition) of the vertex set $V(G)$.  

In 1928 Ramsey \cite{Ram} proved (stated here only for graphs) that for any $r$-coloring of $K_\NN$ there
exists an infinite set $A\subseteq \NN$ such that all of the edges with both endpoints in $A$ have the same color.  It is
natural to wonder if Ramsey's theorem can be strengthened to say something more precise about the size of $A$; for
instance, can we control the growth of the elements of $A$?  It seems to be first noted by Erd\H{o}s in 1964 \cite{E}
that no such result can hold.  That is, there exist $2$-colorings of $K_\NN$ in which every red complete subgraph is finite and the gaps between consecutive vertices of any blue complete subgraph grow arbitrarily fast.

In the case of sparser graphs (such as paths), it may be possible to find monochromatic copies with controlled growth.  In this paper, we prove a variety of such results concerning the \emph{upper density} of a set
$A\subseteq\NN$;
\[
  \bar{d}(A) = \limsup_{n\to\infty} \frac{|A\cap \{1,\ldots,n\}|}{n}.
\]
Given a graph $H$ with vertex set $V(H)\subseteq\NN$, we define the \emph{upper density} of $H$ to be that of $V(H)$.  The \emph{lower density} of $A$, denoted $\ubar{d}(A)$,  is defined analogously in terms of the $\liminf$.  
We note that other more general notions of size, defined via \emph{ideals}\footnote{An \emph{ideal} over a set $X$ is a
non-empty collection $\SI$ of subsets of $X$ which is closed under taking subsets and closed under finite unions.} over $\NN$,
have been studied in the context of Ramsey's theorem going back to the 1970's; the interested reader is directed
to~\cite{Mathias},~\cite{Hrusak}, and~\cite{Farah}.

We note that in order to meaningfully talk about the upper density of a monochromatic copy of an infinite graph $G$ in an edge-colored countably-infinite complete graph $K$, we must fix a ordering of the vertices of $K$ before the coloring of the edges of $K$ is given; as otherwise, for every coloring of the edges of $K$, there exists an ordering of the vertices of $K$ so that $K$ has a complete monochromatic subgraph of density 1.  So for the purposes of computing the upper density, we simply treat the vertex set of $K$ as the positive integers with the natural ordering.

\subsection{Paths}

One of the first results in graph Ramsey theory is the following.
\begin{theorem}[Gerencs\'er, Gy\'arf\'as \cite{GG}]\label{thm:gg}
In every $2$-coloring of $K_n$ there exists a monochromatic path with at least $\ceiling{\frac{2n+1}{3}}$ vertices. Furthermore, this is best possible .
\end{theorem}
Towards an infinite analogue of Theorem~\ref{thm:gg}, one might ask; given a $2$-coloring of $K_\NN$, what is the
largest asymptotic density of a monochromatic path, i.e.\ a monochromatic copy of the graph $P$ with ordered vertex set
$V(P) = \seq{x_n}{n\in\NN}$ and edge set $E(P) = \set{x_n x_{n+1}}{n\in\NN}$?

Erd\H{o}s and Galvin \cite{EG},\cite{EG:R} seem to be the first to consider this question.  Among other
things, they showed that there is a $2$-coloring for which every monochromatic subgraph of bounded maximum degree (in particular, paths) has lower density
$0$ (\cite[Theorem 1.1]{EG}).  This means that we
cannot expect results using the usual asymptotic density.  On the other hand, a result of Rado \cite{R} states that in
any $r$-coloring of $K_\NN$, there is a partition of $\NN$ into at most $r$ monochromatic paths, and this implies that
one of these paths must have upper density at least $1/r$.  In the case of 2 colors, Erd\H{o}s and Galvin strengthened
this lower bound, and also found an upper bound;
\begin{theorem}[Erd\H{o}s, Galvin \cite{EG}]\label{thm:EG2colorUD}~
  In every $2$-coloring of $K_\NN$ there exists a monochromatic path with upper density at least $2/3$.  On the other hand,
  there exists a $2$-coloring of $K_\NN$ in which every monochromatic path has upper density at most $8/9$.\end{theorem}

Erd\H{o}s and Galvin also considered an alternate notion of density of a path $P$ which only measures the density of ``unbroken'' initial segments of $P$.
\begin{definition}
  Given an infinite sequence of distinct natural numbers $A = (a_1,a_2,\ldots, a_n,\ldots )$, the \emph{strong upper density} of $A$ and the \emph{strong lower density}
  of $A$ are defined to be
  \[
    \bar{d}_s(A) = \limsup_{n\to\infty} \frac{f(n)}{n} ~~~\text{ and }~~~
    \underline{d}_s(A) = \liminf_{n\to\infty} \frac{f(n)}{n}
  \]
  respectively, where $f(n) = \sup\set{m}{\{a_1,\ldots,a_m\}\subseteq \{1,\ldots,n\}}$.
\end{definition}

\begin{figure}[ht]
\centering
\includegraphics[scale=0.85]{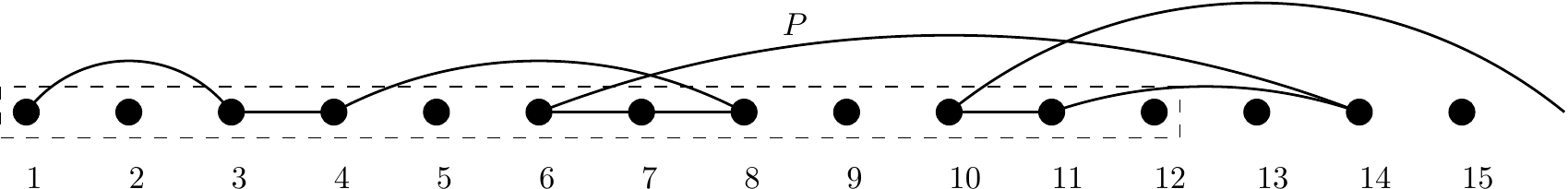}
\caption{A path $P$ which has density $8/12$ in the set $\{1,2,\dots, 12\}$, but which has strong density $6/12$ in the
set $\{1,2,\dots, 12\}$.}\label{fig:sd}
\end{figure}

In other words, in the strong upper density of a sequence $A$ we consider only the density of the longest initial
segment of $A$ which lies entirely inside each set $\{1,\ldots,n\}$ (see Figure~\ref{fig:sd}).  In particular, we have $\bar{d}_s(A)\le
\bar{d}(A)$.  Given a graph $G$ with an ordered vertex set $V(G) = (v_1,v_2,\ldots,v_n,\ldots)$ consisting of distinct
natural numbers, we define the \emph{strong upper density} of $G$ to be the strong upper density of the sequence
$(v_1,v_2,\ldots,v_n,\ldots)$.

Again for 2 colors, Erd\H{o}s and Galvin offered the following bounds on the strong upper density of
monochromatic paths;
\begin{theorem}[Erd\H{o}s, Galvin \cite{EG}]\label{thm:EG2colorSUD}~
  In every $2$-coloring of $K_\NN$ there exists a monochromatic path with strong upper density at least $1/(3 + \sqrt{3})$.  On the other hand,
  there exists a $2$-coloring of $K_\NN$ in which every monochromatic path has strong upper density at most $2/3$.
\end{theorem}

For $r\geq 3$ colors they proved a few results, including the somewhat surprising fact that there exists an $r$-coloring
of $K_\NN$ in which every monochromatic path has strong upper density 0 (\cite[Theorem 1.7]{EG}. 

\begin{theorem}[Erd\H{o}s, Galvin \cite{EG}]  Let $r\geq 3$. 
\begin{enumerate}

\item There exists an $r$-coloring of $K_\NN$ such that for all monochromatic paths $P$, $\bar{d}(P)\leq 2/r$.

\item There exists an $r$-coloring of $K_\NN$ such that for all monochromatic paths $P$, $\bar{d}_s(P)=0$. 
\end{enumerate}
\end{theorem}

We improve the lower bound for ordinary upper density of paths.  

\begin{theorem}\label{mainThm1}
In every $2$-coloring of $K_{\NN}$, there exists a monochromatic path $P$ such that $\bar{d}(P)\geq 3/4$.
\end{theorem}

Furthermore, we prove a tight result for the strong upper density of paths in $2$-colored complete graphs.  As this was the only open case, this completely settles the strong upper density problem for paths.

\begin{theorem}\label{mainThm2}
In every $2$-coloring of $K_{\NN}$, there exists a monochromatic path $P$ such that $\bar{d}_s(P)\geq 2/3$.
\end{theorem}

Finally, we provide an improved upper bound for the upper density of paths when $r\ge 3$ (Corollary \ref{cor:affinepath}).

\subsection{Directed paths}

In an orientation of a path $P=v_1v_2\dots$, we say that $v_i$ is an \emph{out-(in-)switch} if the in-(out-)degree of $v_i$ is 0 in $P$.  Note that any endpoints of $P$ are switches, and that $P$ has one endpoint if $P$ is infinite and two endpoints if $P$ is finite.  We say that $P$ is consistently oriented if the only switches of $P$ are the endpoints, and we say that $P$ is anti-directed if every vertex in $P$ is a switch.  

Raynaud proved a directed analog of Theorem \ref{thm:gg}.

\begin{theorem}[Raynaud \cite{Ray}]
In every $2$-coloring of $\vec{K}_n$, there exists a consistently oriented monochromatic path with at least $n/2+1$ vertices.  Furthermore, this is best possible.  
\end{theorem}

As above with undirected graphs, it is natural to consider an infinite analog of Raynaud's theorem.  So it is somewhat surprising that we prove that for all $\ep>0$, there exists a $2$-coloring of $\vec{K}_\NN$ such that every consistently oriented monochromatic path has upper density less than $\ep$ (Example \ref{ex:directedpath}).  And for $r\geq 3$, there exists an $r$-coloring of $\vec{K}_\NN$ such that every monochromatic consistently oriented path has upper density 0 (Example \ref{ex:directedpath3}).  Note that any density result that holds for consistently oriented paths also holds for paths with finitely many switches (since all of the switches are contained in some finite initial segment of the path).

On the other hand, we show that if $P$ contains infinitely many switches, then in any $r$-coloring of the edges of $\vec{K}_\NN$ (actually, in any $r$-coloring of the edges of any tournament on $\NN$), there exists a monochromatic copy of $P$ with positive upper density (Corollary \ref{cor:infswitch}).  

As mentioned in the previous section, Rado proved that in every $r$-coloring of $K_\NN$, there exists a partition of $\NN$ into at most $r$ monochromatic paths (which implies that some monochromatic path has upper density at least $1/r$).  Actually, his result was more general.  He proved that in every $r$-coloring of $\vec{K}_\NN$, there exists a partition of $\NN$ into at most $r$ monochromatic anti-directed paths (which implies that some monochromatic anti-directed path has upper density at least $1/r$).  So for oriented paths, the bound on the density will depend heavily on the orientation of the path.

\subsection{Connected subgraphs}

An old remark of Erd\H{o}s and Rado says that a graph or its complement is connected, or alternatively, in every $2$-coloring of $K_n$ there exists a monochromatic connected subgraph on $n$ vertices.  This remark was extended by Gy\'arf\'as, and is one of the few examples in graph Ramsey theory for which a tight result is known for infinitely many values of $r$.

\begin{theorem}[Gy\'arf\'as \cite{Gy}]\label{thm:Gy}
Every $r$-coloring of $K_n$ contains a monochromatic connected subgraph on at least $n/(r-1)$ vertices.  Furthermore, this is best possible when $r-1$ is a prime power and $(r-1)^2$ divides $n$.
\end{theorem}

We consider the analogous problem for infinite graphs.  First we need a definition of strong upper density for connected
graphs.

\begin{definition}\label{def:sud}
  Let $T$ be a connected subgraph of $K_{\NN}$.  The \emph{strong upper density} of
  $T$ is
  \[
    \bar{d}_s(T) = \sup_{V\in \pi} \bar{d}_s(V)
  \]
  where $\pi$ is the set of all orderings $(v_1,v_2,\dots)$ of $V(T)$ such that for all $k\geq 1$, the graph $T[\{v_1,
  \dots, v_k\}]$ induced by $\{v_1, \dots, v_k\}$ is connected.
\end{definition}

The reader can easily check that the set $\pi$ in the above definition is
nonempty.  Moreover, we note that, for a connected subgraph $T$ of $K_\NN$, to
show that $\bar{d}_s(T) \geq \alpha$ it suffices to find an increasing sequence of
integers $n_1 < n_2 < \cdots < n_k < \cdots $ such that for each $k$, the
subgraph of $T$ induced by $V(T)\cap \{1,\ldots,n_k\}$ is connected, and
$\limsup_{k\to\infty} |V(T)\cap \{1,\ldots,n_k\}| / n_k \geq \alpha$.

It is not hard to see that the remark of Erd\H{o}s and Rado easily extends to
countably infinite graphs which shows that every 2-coloring of $K_\NN$ contains
a monochromatic connected subgraph $T$ such that $V(T)=\NN$, and thus has
\emph{density} 1; however, Example \ref{ex:affinelower3} shows that upper
density, and not density or lower density, is all we can hope for in the case of
$r\ge 3$ colors.  We show that the infinite analog of Erd\H{o}s and Rado's
remark extends to strong upper density, and Example \ref{ex:stronglower} shows
that we must consider strong upper density rather than ordinary strong density.

\begin{theorem}\label{infTree2}
Every $2$-coloring of $K_\NN$ contains a monochromatic connected graph $T$ such that $\bar{d}_s(T)=1$.
\end{theorem}

A (special case of a) conjecture of Ryser says that in every $r$-coloring of $K_n$, there exists a cover of $V(K_n)$ by at most $r-1$ monochromatic connected subgraphs, and if true this would imply Theorem \ref{thm:Gy}.  Ryser's conjecture is known to be true for $r\leq 5$ (see \cite{Tuza}).  For $r=3$, it is not hard to prove this result and see that it easily extends to countably infinite complete graphs.  This implies that in every 3-coloring of $K_\NN$ there exists a monochromatic connected subgraph with upper density at least $1/2$.  We show that the same holds for \emph{strong} upper density.

\begin{theorem}\label{infTree3}
Every $3$-coloring of $K_\NN$ contains a monochromatic connected graph $T$ such that $\bar{d}_s(T)\geq 1/2$.  
\end{theorem}

As mentioned above, there exists an example to show that Theorem \ref{thm:Gy} is best possible when $r-1$ is a prime power and $(r-1)^2$ divides $n$.  We show that this example extends to infinite graphs (Example \ref{ex:affine}).  

\subsection{0-dense graphs}

For $r\geq 2$, we say that a graph $G$ is \emph{$(0,r)$-dense} if there exists an $r$-coloring of $K_\NN$ such that every monochromatic copy of $G$ has density 0, and if $G$ is $(0,2)$-dense, we just say that $G$ is \emph{$0$-dense}.  If for all $r\geq 2$, $G$ is not $(0,r)$-dense, then we say that $G$ is \emph{positively-dense}.   For directed graphs, we define $(0,r)$-dense and $0$-dense in the analogous way.  
We show that if $G$ has bounded independence number, then $G$ is 0-dense (Proposition \ref{prop:relden}).
We also show that if $G$ has an infinite independent set $B$ and every vertex in $V(G)\setminus B$ has finite degree to $B$, then $G$ is positively-dense (Corollary \ref{cor:AB}).

\section{Notation and outline}
\label{sec:notout}

Given $n\in\NN$ we will write $[n]$ for the set $\{1,\ldots,n\}$.

In many of the results to follow we will make use of a partition of $\NN$ into consecutive intervals $A_1,A_2,\ldots,A_n,\ldots$ such that
\[
  \lim_{n\to\infty} \frac{|A_n|}{|A_1| + \cdots + |A_n|} = 1.
\]
In such cases we will simply call the $A_n$'s a sequence of \emph{fast growing intervals}.  In some situations we will require the sizes $|A_n|$ to grow even more rapidly, according to various conditions, in which case these conditions will be laid out explicitly.  

For a graph $G$ and vertex sets $X,Y\subseteq V(G)$, $e(X,Y)$ represents the
number of edges with one endpoint in $X$ and the other in $Y$, and we let
$e(G)=e(V(G), V(G))$.  We write $\delta(X,Y)$ to mean $\min\{|N(v)\cap Y|: v\in
X\}$.  When $X$ and $Y$ are disjoint sets of vertices we will write $[X,Y]$ for
the complete bipartite graph on parts $X,Y$.

For constants $a$ and $b$, we write $a \ll b$ to mean that given $b$, we can choose $a$ small enough so that $a$ satisfies all of necessary conditions throughout the proof.  More formally, we say that a statement holds for $a\ll b$ if there is a function $f$ such that it holds for every $b$ and every $a\le f(b)$ 
In order to simplify the presentation, we will not determine these functions explicitly. 

Given a finite set $F\subseteq\NN$, we define the \emph{local density} of $F$ to be $d_\ell(F) = |F|/\max(F)$.  We note here that if $A_n$ is a sequence of fast-growing intervals, then $d_\ell(A_n)\to 1$ as $n\to\infty$.

The paper is organized as follows.  We collect all of our examples in Section \ref{sec:example}.  For the convenience of the reader we provide two of the examples of Erd\H{o}s and Galvin, stated here in less general terms to suit our purpose, in Appendix A.  
In Section \ref{sec:udpaths} we prove Theorem \ref{mainThm1} which improves Erd\H{o}s and Galvin's lower bound on the upper density of monochromatic paths in 2-colored $K_\NN$.  
In Section \ref{sec:sudpaths} we prove the more complicated Theorem \ref{mainThm2} which gives a tight bound on the strong upper density of monochromatic paths in 2-colored $K_\NN$.  This result relies on a ``robust'' version of Theorem \ref{thm:gg}, the proof of which uses the ``regularity--blow-up method'' based on Szemer\'edi's regularity lemma \cite{Sz} and Koml\'os, S\'ark\"ozy, and Szemer\'edi's blow-up lemma \cite{KSS}.  
In Section \ref{sec:connected} we prove Theorems \ref{infTree2} and \ref{infTree3} about the strong upper density of monochromatic connected subgraphs in 2- and 3- colored $K_\NN$.  
In Section \ref{sec:0dense} we discuss properties of graphs and digraphs which imply that they are 0-dense or positively dense. 
Finally, in Section \ref{sec:conclusion}, we discuss a large number of open questions stemming from our work.

\section{Examples}\label{sec:example}

\subsection{Directed paths}

\begin{example}\label{ex:directedpath}
For all $\ep>0$, there exists a $2$-coloring of $\vec{K}_\NN$ such that every
monochromatic, consistently-oriented path has upper density less than $\ep$.
\end{example}

\begin{proof}
Let $\ep>0$ and set $k=\ceiling{1/\ep}+1$.  Partition $\NN$ into $k$ infinite
sets $A_1,\dots, A_k$ each of density $1/k$ (for instance, partition based on
residues modulo $k$).  Let $x,y\in \NN$.  If $x,y\in A_i$ and $x<y$, then color
$(x,y)$ red and $(y,x)$ blue.  If $x\in A_i, y\in A_j$ and $i<j$, then color
$(x,y)$ blue and $(y,x)$ red.  Now let $P$ be a monochromatic,
consistently-oriented path.  Then, $P$ has infinite intersection with at most
one of the sets $A_i$; thus $P$ has upper density at most $1/k<\ep$.  
\end{proof}

The natural attempt to strengthen Example~\ref{ex:directedpath}, by partitioning
$\NN$ into infinitely-many infinite sets $A_1,\ldots,A_i,\ldots$ each with
density $0$, fails due to the following fact.
\begin{proposition}
  Suppose $\NN = A_1\cup\cdots\cup A_i\cup\cdots$ is a partition of $\NN$ into
  infinite sets, each with density $0$.  Then there is a set $B$ with density $1$ such
  that $B\cap A_i$ is finite for every $i$.
\end{proposition}
\begin{proof}
  Fix a sequence $\e_i$ tending to zero.
  For each $i$, choose $n_i$ large enough that for all $n\ge n_i$,
  \[
    \frac{1}{n}\card{(A_1\cup \cdots \cup A_i)\cap [n]} < \e_i.
  \]
  Now let
  \[
    B = \bigcup_{i=1}^\infty A_i\cap [n_i]
  \]
  Clearly each $B\cap A_i$ is finite; to see that $B$ has density $1$, note that
  for any given $n$, if $n_i \le n < n_{i+1}$ then $[n]\sm B \subseteq (A_1\cup
  \cdots \cup A_i)\cap [n]$ and hence $B$ has density at least $1 - \e_i$ in
  $[n]$.
\end{proof}

We note here that $3$ colors is enough to ensure not only that
every such path has density $0$, but that such paths grow arbitrarily fast.

\begin{example}\label{ex:directedpath3}
Let $h : \NN\to\NN$, where $h(n)\to\infty$.  Then there exists a $3$-coloring of
$\vec{K}_\NN$ such that for every monochromatic, consistently-oriented path
$P=v_1v_2v_3\dots$, the sequence $(v_1,v_2,\dots)$ eventually grows faster than
$h$.
\end{example}

\begin{proof}
Let $\NN = A_1\cup A_2 \cup \cdots \cup A_n \cup \cdots$ be a partition of $\NN$
into consecutive intervals, where $|A_n|\ge h(n)$.  We color $\vec{K}_\NN$ so
that each $A_n$ is complete in green, and for each $m < n$, the edges from $A_m$
to $A_n$ are red and the edges from $A_n$ to $A_m$ are blue.  Let $P$ be a
monochromatic, consistently-oriented path.  Then for all but finitely-many $n$,
$P$ intersects $A_n$ in at most one vertex.  It follows that for large enough
$k$, the $k$th vertex of $P$ is at least $h(k)$.
\end{proof}

We note that the bounds on the upper density in Examples \ref{ex:directedpath} and \ref{ex:directedpath3} also
apply to any path with only finitely-many switches, as removing any initial segment of the path does not change its upper density.

\subsection{Connected subgraphs}

\begin{example}\label{ex:affine}
Let $r\geq 2$.  If $r-1$ is a prime power, then there exists an $r$-coloring of
$K_\NN$ such that every monochromatic connected subgraph $T$ has upper density
at most $\frac{1}{r-1}$.  
\end{example}

\begin{proof}
Since $r-1$ is a power of a prime, there exists an affine plane of order $r-1$.
Partition $\NN$ into $(r-1)^2$ sets $V_1, \dots, V_{(r-1)^2}$, each of density
$1/(r-1)^2$ (for instance partition based on residues mod $(r-1)^2$).  From the affine plane, we will obtain $r$ partitions (``parallel classes'') $P_1, \dots, P_r$ of $\{V_1, \dots, V_{(r-1)^2}\}$ such that each partition consists of $r-1$ classes (``lines'') each containing $r-1$ elements (``points'') and for all $1\leq i< j\leq (r-1)^2$ there exists exactly one $k$ such that $V_i$ and $V_j$ are in the same partition class in $P_k$.  Now for all $1\leq i<j\leq (r-1)^2$ we color all edges between $V_i$ and $V_j$ with color $k$ if and only if $V_i$ and $V_j$ are in the same partition class in $P_k$, and for all $1\leq i\leq (r-1)^2$, we color the edges inside the set $V_i$ arbitrarily.  

So any monochromatic connected subgraph intersects at most $r-1$ of the sets
$V_1, \dots, V_{(r-1)^2}$ and thus has density at most $(r-1)/(r-1)^2=1/(r-1)$.  
\end{proof}

\begin{corollary}\label{cor:affinepath}
Let $r\geq 3$ and let $q$ be the largest prime power which is at most $r-1$.
There exists an $r$-coloring of $K_\NN$ such that every monochromatic path $P$
has upper density at most $1/q<2/r$.   
\end{corollary}

\begin{proof}
To verify that $1/q<2/r$ for $r\geq 3$, first note that this holds for $r=3$ and then let $r\geq 4$.  By Bertrand's postulate\footnote{Using stronger results than Bertrand's postulate, it is possible to get a better upper bound on $1/q$ for sufficiently large $r$.}, there exists a prime $q'$ such that $\floor{r/2}<q'<2\floor{r/2}$.  Choosing $q$ to be the largest such prime, we have $r/2<q\leq r-1$ and thus $1/q<2/r$.  Now take the $(q+1)$-coloring of $K_\NN$ given in Example \ref{ex:affine}, which is an $r$-coloring since $q+1\leq r$.  
\end{proof}

The next two examples justify why we only consider (strong) upper density for monochromatic connected subgraphs.

\begin{example}\label{ex:stronglower}
There exists a $2$-coloring of $K_\NN$ such that every monochromatic connected
subgraph $T$ has strong lower density $0$.  
\end{example}

\begin{proof}
Partition $\NN$ into fast growing intervals $A_1, A_2, \dots$.  For all $x\in A_{2i+1}$ color all edges from $x$ to $y\in A_j$ with $j\leq 2i+1$ red.   For all $x\in A_{2i}$ color all edges from $x$ to $y\in A_j$ with $j\leq 2i$ blue.  If $R$ is a red connected subgraph and $B$ is a blue connected subgraph then 
\[
\ubar{d}_s(R)=\lim_{i\to \infty}\frac{|A_1\cup A_3\cup \dots \cup A_{2i-1}|}{|A_1\cup A_2\cup \dots\cup A_{2i}|}=0 ~~\text{ and }~~ \ubar{d}_s(B)=\lim_{i\to \infty}\frac{|A_2\cup A_4\cup \dots \cup A_{2i}|}{|A_1\cup A_2 \dots\cup A_{2i+1}|}=0
\]
\end{proof}

\begin{example}\label{ex:affinelower3}
There exists a $3$-coloring of $K_\NN$ such that every monochromatic connected
subgraph $T$ has lower density $0$.  
\end{example}

\begin{proof}
Partition $\NN$ into fast growing intervals $A_0, A_1, A_2, \ldots$ and for each
$q = 0,1,2$ and $3$ let $B_q$ be the union of all $A_n$ such that $n\equiv q\pmod{4}$.
Note that each $B_q$ has upper density $1$.  Let all edges between $B_0$ and
$B_3$ and between $B_1$ and $B_2$ be red, all edges between $B_0$ and $B_2$ and
between $B_1$ and $B_3$ be blue, and all edges between $B_0$ and $B_1$ and between
$B_2$ and $B_3$ be green.  Color all edges inside each $B_q$ arbitrarily.  Then any
monochromatic, connected subgraph $T$ is disjoint from some $B_q$; for instance
any red connected subgraph must either be contained in $B_0\cup B_3$ or $B_1\cup
B_2$.  Hence any monochromatic connected subgraph has lower density $0$.
\end{proof}

\section{Upper density of paths}\label{sec:udpaths}

In this section we will prove Theorem \ref{mainThm1}; that is, in every $2$-coloring of $K_{\NN}$, there exists a monochromatic path $P$ such that $\bar{d}(P)\geq 3/4$.

We will use the following result implicitly proved by Gy\'arf\'as and Lehel (see Pokrovskiy \cite[Corollary 1.9]{P}).

\begin{theorem}[Gy\'arf\'as, Lehel \cite{GL}]\label{thm:bipaths}
In every $2$-coloring of $K_{n,n}$ there exists a partition of the vertex set into at most three disjoint monochromatic paths.
\end{theorem}

Let $G$ be a graph.  A \emph{total $r$-coloring of $G$} is a coloring of both
the edges and vertices of $G$ with a set of $r$ colors.  Given a total
$r$-coloring of $G$, a \emph{monochromatic path forest} of color $c$ is a collection of
vertex disjoint paths (i.e.\ a path forest) such that every edge and the
endpoints of every path have color $c$.  Theorem \ref{thm:bipaths} implies
the following result about monochromatic path forests in complete
graphs with a total $2$-coloring.

\begin{lemma}
  \label{lemma:GLP}
  Let a total $2$-coloring of $K_n$ be given where, if $R$ is the set of red vertices and $B$ the set of blue
  vertices, then $|R| \le |B|$.  Then for every $B'\subseteq B$ with $|B'|=|R|$, there exist red and blue monochromatic path forests $F_R$ and $F_B$ such that $F_R\subseteq R\cup B'$ and 
  \[
    |F_R| + |F_B| \ge n + |R| - 3.
  \]
\end{lemma}

\begin{proof}
Let $B'\subseteq B$ with $|B'|=|R|$ and apply Theorem \ref{thm:bipaths} to the bipartite graph induced by $B'$ and $R$ to obtain a partition of $B'\cup R$ into $1\leq t\leq 3$ monochromatic paths $P_1, \dots, P_t$.  For all $i\in [t]$, if $P_i$ is red, let $d_i$ be the number of endpoints of $P_i$ which are in $B'$, and if 
$P_i$ is blue, let $d_i$ be the number of endpoints of $P_i$ which are in $R$.  Let $d:=d_1+\dots+d_t$.  Now suppose we have chosen the paths $P_1, \dots, P_t$ satisfying Theorem \ref{thm:bipaths} in such a way that $d$ is minimum.  This implies that we can't have a red path $P_i$ with $d_i>0$ and a blue path $P_j$ with $d_j>0$, because by using the edge between their endpoints, we can obtain a partition into monochromatic paths with a smaller value of $d$.  

If there exists $P_i$, say blue, with $d_i=2$, then since $[B', R]$ is balanced, none of the paths $P_1, \dots, P_t$ are red and $d\leq 3$.  Otherwise $d_i\leq 1$ for all $i\in [t]$ and thus $d\leq 3$.   So by deleting the at most 3 vertices which are endpoints in the wrong set and letting $F_R$ be the vertices in $R$ together with the red paths and letting $F_B$ be the vertices in $B$
together with the blue paths gives $$|F_R|+|F_B|\geq |B|+|R|+|B'|-3=n+|R|-3,$$ as desired.
\end{proof}

Now we will use the Lemma above with a ``density increment''-type argument to get a monochromatic path forest with density arbitrarily close to $3/4$.  

\begin{lemma}
  \label{lemma:mpf}
  Let $0<\ep\leq 3/4$ and $k\in\NN$ with $k\geq 3/\ep$.  For all $n\ge 4k/\ep$ and for every total $2$-coloring of the
  complete graph $K_n$ on the vertex set $[n]$, there exists a monochromatic path forest $F$ such that
  for some $\ell\ge k$, $\card{F\cap [\ell]} \ge (3/4 - \e)\ell$.
\end{lemma}

\begin{proof}
  Fix a total $2$-coloring of the complete graph on the vertex set $[n]$, and let $R$ and $B$ be
  the sets of red vertices and blue vertices respectively.  We will assume without loss of generality that $|R| \le
  |B|$.  

  We may assume that for each $\ell\ge k$, $|R\cap [\ell]|$ is less than $|B\cap [\ell]|$; otherwise there is some
  $\ell\ge k$ such that $|R\cap [\ell]| = |B\cap [\ell]|$, and then by Lemma~\ref{lemma:GLP} there is a
  monochromatic path forest $F\subseteq [\ell]$ such that
  \[
    |F| \ge \frac{1}{2}\left(\ell+\frac{\ell}{2}-3\right)=\left(\frac{3}{4}-\frac{3}{2\ell}\right)\ell\geq
    \left(\frac{3}{4}-\ep\right)\ell.
  \]

  Let $b_0=|B|$, $r_0 = |R|$, and for all $i$ such that $r_i > 0$, let $b_{i+1}=r_i$ and let $r_{i+1}$ be the number of red
  vertices before the $b_{i+1}$th blue vertex, according to the natural order of $[n]$.    For each $i\ge 0$, let $J_i =
  [r_i+b_i]$.  Note that $J_i$ has exactly $r_i$ red vertices and exactly $b_i$ blue vertices; hence, by the assumption
  made above, for each $i$ we have $r_{i+1} < b_{i+1} = r_i$.  It follows that the sizes of the intervals $J_i$ are
  decreasing.  Let $t$ be the largest integer such that $|J_t| \ge k$.  Note that by the maximality of $t$, we must have $r_t<k$,
  as otherwise $b_{t+1}=r_t\geq k$ and thus $r_{t+1}+b_{t+1}\geq k$.  Moreover, if $b_t \ge 3 r_t$, then the blue
  vertices in $J_t$ form a monochromatic path forest with density at least $3/4$ in $J_t$; hence we may assume that 
  \begin{equation}\label{btrt}
    b_t - r_t < 2 r_t < 2k.
  \end{equation}
  
  The purpose of the following claim is to show that, if there is no monochromatic path forest with the desired density
  in some interval $J_i$, then the differences $b_i - r_i$ increase as $i$ increases.  Combined with \eqref{btrt}, it
  will follow that $b_0 - r_0$ is very small relative to $n$, and we will get a sufficiently dense monochromatic path
  forest by applying Lemma~\ref{lemma:GLP} again.
  
  \begin{claim}\label{claim:ribi}
    For each $i$ with $1\le i\le t$, either
    \begin{enumerate}
      \item  there is a blue monochromatic path forest $F_B$ with $\card{F_B\cap J_{i-1}} \ge (3/4 - \e)|J_{i-1}|$,
      \item  there is a red monochromatic path forest $F_R$ with $\card{F_R\cap J_i} \ge (3/4 - \e)|J_i|$, or
      \item  $b_{i-1} - r_{i-1} < b_i - r_i$.
    \end{enumerate}
  \end{claim}
  
  \begin{proof}
    Let $1\leq i\leq t$ and consider the interval $J_{i-1}$.  By Lemma~\ref{lemma:GLP}, there are red and blue
    monochromatic path forests $F_R$ and $F_B$ contained in $J_{i-1}$ such that $|F_R| + |F_B| = b_{i-1} + 2r_{i-1} - 3$, and
    $F_R\cap B$ is contained in the first $b_i = r_{i-1}$-many blue vertices of $J_{i-1}$.  If $|F_B| \ge (3/4 -
    \e)|J_{i-1}|$ then we are done.  Assuming otherwise, we have
    \begin{align*}
      |F_R\cap B|  =  b_{i-1} + 2r_{i-1} - |F_B| - |F_R\cap R| - 3 & = b_{i-1} + r_{i-1} - |F_B| - 3 \\
        & > \left(\frac{1}{4}+\ep\right)(b_{i-1} + r_{i-1}) - 3 \\
        & \ge \frac{1}{4} (b_{i-1} + r_{i-1}),
    \end{align*}
    where the last inequality holds by the definition of $k$ and the choice of $t$.  
    Now consider the size of $F_R\cap J_i$.  Since $F_R\cap B\subseteq J_i$, we have
    \[
      |F_R\cap J_i| > r_i + \frac{1}{4}(b_{i-1} + r_{i-1}).
    \]
    If $|F_R\cap J_i| \ge (3/4 - \e)|J_i|$, then we are done.  Assuming otherwise, then, we get
    \[
      r_i + \frac{1}{4}(b_{i-1} + r_{i-1}) < \frac{3}{4}(b_i + r_i)=\frac{1}{2}r_{i-1}+\frac{1}{4}b_i+\frac{3}{4}r_i
    \]
    Rearranging, we get the desired inequality $b_{i-1} - r_{i-1} < b_i - r_i$.
   \end{proof}
  
  By Claim~\ref{claim:ribi}, we may assume that $b_{i-1} - r_{i-1} < b_i - r_i$ for each $i$ with $1\le i\le t$, and
  thus by \eqref{btrt}, we have 
  $b_0 - r_0 < b_t - r_t<2k$.  So by Lemma~\ref{lemma:GLP}, we may find a monochromatic path forest $F\subseteq J_0
  = [n]$ with size at least
  \[
    \frac{1}{2}(|B| + |R| + |R| - 3)\ge \frac{1}{2}(3|B| - 4k - 3)\ge \frac{1}{2}\left(\frac{3n}{2} - \e n - 3\right) >
    \left(\frac{3}{4} - \e\right) n.
  \]
  This completes the proof.
\end{proof}

\begin{proof}[Proof of Theorem \ref{mainThm1}]
  Fix a $2$-coloring of the edges of $K_\NN$.  We define a $2$-coloring of the vertices by coloring $n\in\NN$ red if and
  only if there are infinitely-many $m\in\NN$ such that the edge $\{n,m\}$ is red.

  Case 1:\footnote{The proof in this case is exactly as in~\cite[Theorem~3.5]{EG}.} there are vertices $x$ and $y$ of the same color,
  say red, and a finite set $S\subseteq\NN$ such that there is no red path disjoint from $S$ which connects $x$ to $y$.
  We partition $\NN\setminus S$ into sets $X, Y, Z$, where $x'\in X$ if and only if there is a red path, disjoint from
  $S$, which connects $x$ to $x'$, and $y'\in Y$ if and only if there is a red path disjoint from $S$ which connects $y$
  to $y'$.  Note that every edge from $X\cup Y$ to $Z$ is blue.  Since $x$ and $y$ are colored red, both $X$ and $Y$ are infinite, and by choice of $x$ and $y$ all edges in
  the bipartite graph between $X$ and $Y\cup Z$ are blue.  Hence there is a blue path with vertex set $X\cup Y\cup Z = \NN\setminus
  S$.

  Case 2: for every pair of vertices $x$ and $y$ of the same color $\chi$, and every finite set $S\subseteq\NN$, there is
  a path from $x$ to $y$ of color $\chi$ which is disjoint from $S$.  Let $\e_n$ be a sequence of positive reals tending
  to zero, and let $a_n$ and $k_n$ ($n\in\NN$) be increasing
  sequences of integers such that
  \begin{itemize}
    \item  $k_n \ge 3/\e_n$ and $a_n \ge 4 k_n/\e_n$, and
    \item  as $n\to \infty$, $k_n / (a_1 + \cdots + a_{n-1} + k_n) \to 1$.
  \end{itemize}
  Let $\NN = A_1\cup A_2 \cup \cdots \cup A_n \cup \cdots$ be a partition of $\NN$ into consecutive intervals such that
  $|A_n| = a_n$.  By Lemma~\ref{lemma:mpf} we may find monochromatic path forests $F_n \subseteq A_n$ and initial segments
  $I_n\subseteq A_n$ of length $\ge k_n$ such that 
  \[
    \card{F_n\cap I_n} \ge \left(\frac{3}{4} - \e_n\right) \card{I_n}.
  \]
  It follows that for any $B\subseteq\NN$ containing infinitely-many $F_n$'s, we have
  \[
    \bar{d}(B) \ge \limsup_{n\to\infty} \frac{\card{F_n\cap I_n}}{a_1 + \cdots + a_{n-1} + |I_n|} = 
    \limsup_{n\to\infty} \frac{3}{4} - \e_n = \frac{3}{4}
  \]
  By the pigeonhole principle, there are infinitely-many $F_n$'s of the same color, say blue.  We will recursively
  construct a blue path $P$ which contains infinitely-many of these $F_n$'s.  To
  see how this is done, suppose we have
  constructed a finite initial segment $p$ of $P$.  We will assume as an inductive hypothesis that $p$ ends on one of
  the end-vertices of some blue $F_m$.  Let $v$ be this end-vertex.  Let $n$ be large enough that $\min{A_n}$ is greater
  than every vertex in $p$, and $F_n$ is blue.  Let $w_1,\ldots,w_s$ list the end-vertices of $F_n$.  By the case
  assumption, there is a blue path $q_1$ connecting $v$ to $w_1$, such that $q_1$ is disjoint from $A_1\cup \cdots \cup
  A_n$.  Similarly, there is a blue path $q_2$ connecting $w_1$ to $w_2$, such that $q_2$ is disjoint from $A_1\cup
  \cdots \cup A_n \cup q_1$.  Continue in this fashion, choosing disjoint blue paths $q_3,\ldots,q_s$ such that $q_i$
  connects $w_{i-1}$ to $w_i$.  Clearly the paths of $F_n$ allow us to put all of these together into a blue path $p'$
  which extends $p$ and contains all of the vertices of $F_n$.
\end{proof}

\begin{remark}
  One may use an ultrafilter in the proof of Theorem~\ref{mainThm1} to define the $2$-coloring of $\NN$ in such a way
  that Case 1 never occurs.  (See Theorem~\ref{thm:ultratournament} for an example of this type of argument in action.)
 While this shortens the proof somewhat, the above proof has the advantage of being more elementary.
\end{remark}

\section{Strong upper density of paths}\label{sec:sudpaths}
In this section, we will prove Theorem \ref{mainThm2}; that is, in every $2$-coloring of $K_{\NN}$, there exists a monochromatic path $P$ such that $\bar{d}_s(P)\geq 2/3$.

To prepare for the proof of the main result in Section \ref{mainsubsection}, we first develop a number of tools in the next two subsections.

\subsection{General Lemmas}

Given a graph $G$ and disjoint sets $X,Y\subseteq V(G)$, define $d(X,Y)=\frac{e(X,Y)}{|X||Y|}$.  Given $\ep\geq 0$, say that a bipartite graph with parts $X,Y$ is $\ep$-\emph{regular} if $|d(X',Y')-d(X,Y)|\leq \ep$ for all $X'\subseteq X$ and $Y'\subseteq Y$ with $|X'|> \ep |X|$ and $|Y'|>\ep |Y|$.  Given $\ep,d\geq 0$, say that a bipartite graph $H$ with parts $X,Y$ is $(\ep,d)$-\emph{super-regular} if $H$ is $\ep$-regular and $\delta(X,Y)\geq d|Y|$ and $\delta(Y,X)\geq d|X|$.  An $(\ep, d, m)$-\emph{super-regular cycle} of order $k$ is a graph obtained from taking a cycle $C$ of order $k$ (allowing for an edge to be a cycle of order $2$) and replacing each vertex of $C$ with a set of $m$ vertices and replacing each edge of $C$ with an $(\ep, d)$-super-regular bipartite graph.  

Erd\H{o}s and Galvin proved in \cite{EG} that every 2-coloring of $K_\NN$
contains a monochromatic path of strong upper density at least $1/(3+\sqrt{3})$.
One of the key elements of their proof is the use of Gerencs\'er and
Gy\'arf\'as's theorem (Theorem \ref{thm:gg}) that every $2$-coloring of $K_n$
admits a monochromatic path with at least $2n/3$ vertices.  One of the key
elements of our proof of Theorem \ref{mainThm2} is the following ``robust''
version of Theorem \ref{thm:gg}.

\begin{lemma}\label{supercycles}
For all $0<\eta\leq 1/4$ there exists $0<\frac{1}{n_0}\ll \frac{1}{m}\ll \ep< 4d\ll\sigma$ such that every $2$-colored graph $K_n$ on $n\geq n_0$ vertices has a monochromatic $(\ep, d, m)$-super-regular cycle $C_1$ on at least $(2/3-\eta)n$ vertices.
\end{lemma}

The proof of Lemma \ref{supercycles} is a straightforward consequence of the regularity/blow-up method combined with the following result of Gy\'arf\'as and S\'ark\"ozy.\footnote{We note that the statement of the above result in \cite{GS} is actually for paths instead of cycles; however, their proof implicitly gives the result for cycles.  Furthermore, Benevides, \L uczak, Skokan, Scott, and White \cite{BLSSW} proved a stronger version of Theorem \ref{conmatch} which says that if $\frac{1}{n_0}\ll \eta\leq 1/4$ and $G$ is a graph on $n\geq n_0$ vertices with $\delta(G)\ge 3n/4$, then in every $2$-coloring of $G$, there exists a monochromatic cycle of length at least $(2/3-\eta)n$.}

\begin{theorem}[Gy\'arf\'as, S\'ark\"ozy \cite{GS}]\label{conmatch}
Let $\frac{1}{n_0}\ll \eta\leq 1/4$.  If $\delta(G)\geq (3/4+\eta)n$, then in every $2$-coloring of $G$, there exists a monochromatic cycle of length at least $(2/3-\eta)n$.
\end{theorem}

A proof of Lemma \ref{supercycles} can be found in Appendix B, along with a proof of the following consequence of the  well-known Blow-up lemma of Koml\'os, S\'ark\"ozy, and Szemer\'edi \cite{KSS}.

\begin{lemma}\label{blowupcor}
Let $k\geq 2$ and $0<\frac{1}{m}\ll \ep, d$.  An $(\ep, d, m)$-super-regular cycle $\mathcal{C}$ of order $k$ has the property that for all $u,v\in V(\mathcal{C})$, there exists a $u,v$-path of length at least $(m-1)k$.
\end{lemma}

\subsection{Specific Lemmas}

\begin{definition}
\label{definition:maximal-robust-cycle}
Let $G$ be a graph on $n$ vertices and let $\alpha \in [0,1]$.  We say that $X\subseteq G$ is an
\emph{$\alpha$-connector} if, for all distinct $u,v\in V(X)$, there is a $u,v$-path $P\subseteq X$, such that $|V(P)|\ge
\alpha n$.

We say that that an $\alpha$-connector $X$ is \emph{maximal} if, for all $u\not\in V(X)$, $u$ has at most one neighbor in $X$.
\end{definition}

We now show how to use a $(\ep, d, m)$-super-regular cycle to construct a maximal $\alpha$-connector.

\begin{lemma}
  \label{lemma:maximal-robust-cycles}
  Let $0<\frac{1}{n_0}\ll \ep\leq 2/3$.  For all 2-colored complete graphs on $n\geq n_0$ vertices,  there exists a monochromatic maximal $\alpha$-connector $X$ with $\alpha \ge 2/3 - \e$.
\end{lemma}

\begin{proof}
Apply Lemma \ref{supercycles} to get a monochromatic, say blue, super-regular cycle $X_0$ on at least $(2/3-\ep/2)n$ vertices.  Now for each $i$, let $X_{i+1}$ be the set of vertices which have at least two blue edges into $\bigcup_{j=0}^i X_j$.  Note that $X_i\subseteq
X_{i+1}$.  Now let $X = \bigcup_{i=0}^\infty X_i$. 

If $X$ is an $\alpha$-connector, then by definition it will be maximal.  So to see that $X$ is an $\alpha$-connector, let $u_i, v_j\in X$ where $u_i\in X_i$ and $v_j\in X_j$.  By the
definition of $X_1, X_2, \dots$, it is possible to iteratively choose disjoint paths from $u_{i}$ to some vertex $u_0\in X_0$ and
from $v_j$ to some vertex $v_0\in X_0$ with $u_0\neq v_0$.  Now using the fact that $X_0$ is a super-regular cycle on at least $(2/3-\ep/2)n$ vertices, we use Lemma \ref{blowupcor} find a blue path from $u_0$ to $v_0$ which uses nearly all of the vertices of $X_0$, i.e.\ a blue $u_i, v_j$-path on at least $(2/3-\ep)n$ vertices.
\end{proof}

In the next lemma, we make use the of the following result to show that a nearly complete bipartite graph has nearly spanning paths with endpoints in a prescribed set.

\begin{theorem}[Las Vergnas \cite{LV} (see Berge {\cite[Chapter 10, Theorem 14]{Ber}})]\label{berge_lemma}
Let $G=(U,V,E)$ be a bipartite graph on $2m\geq 4$ vertices with vertices in $U=\{u_1, \dots, u_m\}$ and $V=\{v_1, \dots, v_m\}$ such that $d(u_1) \leq \cdots \leq d(u_m)$ and $d(v_1) \leq \cdots \leq d(v_m)$.  If for the smallest two indices $j$ and $k$ such that $d(u_j)\leq j+1$ and $d(v_k)\leq k+1$ (if they exist), we have $$d(u_j)+d(v_k)\geq m+2,$$
then for all $u\in U$ and $v\in V$, there exists a Hamiltonian path having $u$ and $v$ as endpoints.
\end{theorem}

\begin{lemma}\label{lem:allinone}
Let $0<\frac{1}{n_0}\ll \ep\ll \alpha_1, \alpha_2\leq 1$ where $\alpha_2 > 1/2 +
\e$ and let $G$ be a 2-colored complete graph on a finite interval
$V\subseteq\NN$ with $|V|\geq n_0$.  Let $V_1\cup V_2$ be a partition of $V$
such that $V_1$ is an initial segment of $V$, and both $V_1$ and $V_2$ have
local density\footnote{See Section~\ref{sec:notout}} at least $1-\ep$.  Suppose
$X_2$ is a maximal blue $\alpha_2$-connector in $V_2$.
\begin{enumerate}
\item If $X_1\subseteq V_1$ is such that $|X_1|/|V_1|\geq \alpha_1$ and there is no blue matching of size 2 from
$X_1$ to $X_2$, then for each $i\in [2]$, there is a red path $P\subseteq V_1\cup V_2$ with both endpoints in $X_i$ such that $P$ has local
density at least $(1-\ep)^2\alpha_1$.

\item If $X_1$ is a maximal red $\alpha_1$-connector in $V_1$, then there exist $u_i,v_i\in X_i$ for some $i\in [2]$ and
red and blue $u_i,v_i$-paths $P_1$ and $P_2$ respectively, such that $P_j$ has local density at least $(1 -
\ep)^2\alpha_j$ for each $j\in [2]$.
\end{enumerate}
\end{lemma}

\begin{proof}
(i) Note that since $V_2$ has local density at least $1 - \ep$, we have $|V_1|
\le \ep |V_2|$.  Since $|X_2| \ge \alpha_2|V_2|$ and $\alpha_2 > 1/2 + \e$, it
follows that $|X_2| > |X_1| + |V_2\setminus X_2|$.  We may thus find an initial
segment $I$ of $V_2$ such that $|X_1| + |(V_2\setminus X_2)\cap I| = |X_2\cap
I|$, by beginning with $I = \emptyset$ and adding the elements of $V_2$ in order
until this equality is achieved.  Let $Y_1 = X_1\cup ((V_2\setminus X_2)\cap I)$
and $Y_2 = X_2\cap I$.  By the assumption, there is no blue matching of size 2 from $X_1$ to $X_2$, and by the fact that $X_2$ is a maximal
$\alpha_2$-connector in $V_2$, every vertex in $(V_2\setminus X_2)\cap I$ has at most 1 blue neighbor in $X_2$.  So there exists a vertex $y_1'\in Y_1$ and a vertex $y_2'\in Y_2$ such that the bipartite graph $[Y_1\setminus\{y_1'\}, Y_2\setminus \{y_2'\}]$ has at least $(|Y_1|-1)(|Y_2|-1)-(|Y_1|-|X_1|)$ red edges.  It follows from Theorem \ref{berge_lemma} (for instance) that for
each $i$ we can find a red path $P_i$ with both endpoints in $Y_i$ which misses
at most 3 vertices from $Y_1\cup Y_2 = X_1\cup I$.  Then $P_i$ has local density
at least
\[
  \frac{|X_1|+|I|-3}{\min(V) + |V_1| + |I|}\ge \frac{|X_1| - 3}{\min(V) + |V_1|}
  \ge (1 - \e)\frac{|X_1| - 3}{|V_1|} \ge (1 - \e)^2\alpha_1.
\]

(ii) 
Case 1: There is a blue matching $\{u_1 u_2, v_1 v_2\}$ of size $2$ between $X_1$ and $X_2$.  We construct $P_1$ and
$P_2$ as follows.  Since $X_1$ is a red $\alpha_1$-connector inside $G[V_1]$, we may choose a red path $P_1$ with
endpoints $u_1$ and $v_1$ such that $\card{P_1}\ge \alpha_1 \card{V_1}$.  Similarly, we may choose a blue path $P_2^*
\subseteq X_2$ with endpoints $u_2$ and $v_2$, such that $\card{P_2^*} \ge \alpha_2\card{V_2}$.  It follows that
$d_\ell(P_1)\ge (1-\ep)\alpha_1$ and $d_\ell(P_2^*)\ge (1 - \ep)\alpha_2$.  Now we let $u = u_1$, $v = v_1$, and $P_2 =
u_1 u_2 P_2^* v_2 v_1$.

Case 2: There is no blue matching of size $2$ between $X_1$ and $X_2$.  Then we apply  Lemma \ref{lem:allinone}.(i) to
get a red path $P_1$ with endpoints $u,v\in X_2$ such that $P_1$ has local density at least $(1-\ep)^2\alpha_1$.  Since
$X_2$ is a blue $\alpha_2$-connector, we may also choose a blue path $P_2$ with endpoints $u$ and $v$ and
$\card{P_2}\ge\alpha_2\card{V_2}$, and it follows that the local density of $P_2$ is at least $(1 - \ep)\alpha_2$.  
\end{proof}

\subsection{Main result}\label{mainsubsection}

\begin{proof}[Proof of Theorem \ref{mainThm2}]
  Fix some sequence $\e_n > 0$ tending to zero.  By
  Lemma~\ref{lemma:maximal-robust-cycles}, there is a sequence $k_n$ of integers
  such that every 2-colored complete graph on at least $k_n$ vertices admits a
  monochromatic maximal $\alpha$-connector for some $\alpha\ge 2/3 - \e_n$.
  Partition $\NN$ into a sequence $A_n$ of fast-growing intervals such that
  $|A_n|\ge k_n$ and moreover the local density of $A_n$ is at least $1 - \e_n$.
  Now fix a 2-coloring of $K_\NN$.  Then, for each $n$ we may find a maximal
  $\alpha_n$-connector $X_n$ in $A_n$ of some color $\chi_n$, with $\alpha_n\ge
  2/3 - \e_n$.  (See Figure \ref{fig:cases} in each of the cases below)

\begin{figure}[ht]
    \centering
    \begin{subfigure}[b]{0.35\textwidth}
        \includegraphics[width=\textwidth]{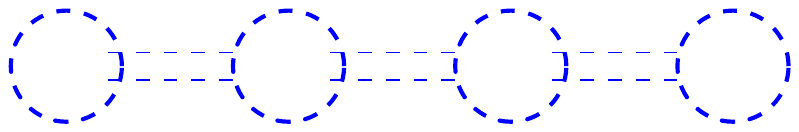}
        \caption{Case 1a}
    \end{subfigure}
    \qquad
    \begin{subfigure}[b]{0.55\textwidth}
        \includegraphics[width=\textwidth]{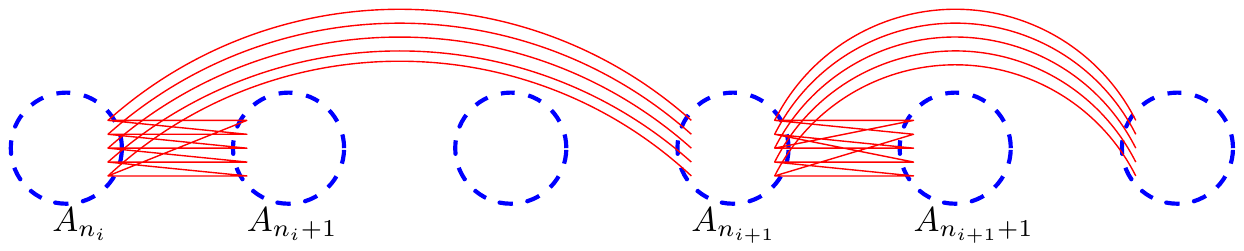}
        \caption{Case 1b}
    \end{subfigure}
    \\
     \vspace{.2cm}
    \begin{subfigure}[b]{0.9\textwidth}
        \includegraphics[width=\textwidth]{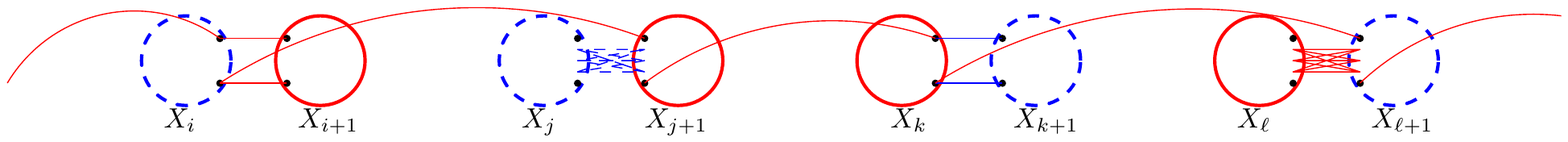}
        \caption{Case 2}
    \end{subfigure}
 
    \caption{Cases in the proof of Theorem \ref{mainThm2}}\label{fig:cases}
\end{figure}

  Case 1: $\chi_n$ is eventually constant.  Say, without loss of generality, that $\chi_n$ is eventually blue.  
  
  Case 1(a): There is an infinite sequence $n_i$ in $\NN$ such that for each $i$, there is a blue matching of size $2$
  between vertices in $X_{n_i}$ and $X_{n_{i+1}}$.  Using these matchings we may choose distinct vertices $u_i,v_i\in
  X_{n_i}$ such that for each $i$, the edge $v_i u_{i+1}$ is blue.  By Definition~\ref{definition:maximal-robust-cycle},
  we may choose blue paths $P_i \subseteq X_{n_i}$, with endpoints $u_i$ and $v_i$, such that $\card{P_i} \ge (2/3 -
  \e_{n_i}) \card{A_{n_i}}$, and hence the local density of $P_i$ is at least $(2/3 - \e_{n_i})(1 - \e_{n_i})$. Then by
  joining the paths $P_i$ with the edges $v_i u_{i+1}$, we get a blue path $P$ with $\bar{d}_s(P) \ge 2/3$.

  Case 1(b): There is an infinite sequence $n_i$ in $\NN$ such that for each $i$ and for each $m > n_i$, there is no
  blue matching of size $2$ between vertices in $X_{n_i}$ and $X_m$.  We may assume that for each $i$, $n_{i+1} > n_i +
  1$.  Note that the induced graph on $A_{n_i}\cup A_{n_i+1}$ satisfies the hypotheses of Lemma \ref{lem:allinone}.(i).
  So there exists a red path $P_i$ having both endpoints $u_i,v_i$ in $A_{n_i}$ such that the local density of $P_i$ is at least
  $(1-\ep_{n_i})^3(2/3-\ep_{n_i}).$ Now by joining the $P_i$'s with the edges $v_i u_{i+1}$, we get a red path $P$ with
  $\bar{d}_s(P)\ge 2/3$.

  Case 2: $\chi_n$ is not eventually constant.  Then there is an infinite set $M\subseteq\NN$ such that for all $n\in
  M$, $\chi_n \neq \chi_{n+1}$.  By Lemma~\ref{lem:allinone}.(ii), for each $n\in M$ we may choose vertices
  $u_n,v_n\in X_n\cup X_{n+1}$ and monochromatic paths $P_n$ and $Q_n$ such that
  \begin{enumerate}
    \item  $P_n$ and $Q_n$ both have endpoints $u_n,v_n$,
    \item  $P_n$ is blue and $Q_n$ is red, and
    \item  each of $P_n$, $Q_n$ has local density at least
    \[
      (1 - \e_n)^2(2/3 - \e_n).
    \]
  \end{enumerate}
  By Ramsey's theorem, we may find an infinite sequence $n_i$ in $M$ such that for each $i$, the edge $v_{n_i}u_{n_{i+1}}$ is the same color $\rho$.  Using a combination of these edges and either the $P_{n_i}$'s or the
  $Q_{n_i}$'s, we may construct an infinite path $P$ with color $\rho$, and
  \[
    \bar{d}_s(P) \ge \limsup_{i\to\infty} (1 - \e_{n_i})^2(2/3 - \e_{n_i}) = 2/3,
  \]
  which completes the proof.
\end{proof}

\section{Strong upper density of connected subgraphs}\label{sec:connected}

The main result of this section is Theorem \ref{infTree3}; that is, in every 3-coloring of the edges of $K_{\NN}$ there exists a monochromatic connected graph $T$ such that $\bar{d}_s(T)\geq 1/2$.  

Before getting to the main result, we prove the following lemma.  Note that a  \emph{maximal monochromatic connected subgraph} is a monochromatic connected subgraph whose vertex set is not properly contained in the vertex set of a monochromatic connected subgraph of any color.

\begin{lemma}
  \label{lem:3col}
Let $K$ be a complete graph (on a finite or infinite vertex set). For every 3-coloring of $K$, either
\begin{enumerate}
\item there exists a monochromatic connected subgraph on $n$ vertices, or

\item there exists a partition $\{W, X, Y, Z\}$ of $[n]$ (all parts non-empty),
such that $B^1:=[W,X]$ and $B^2:=[Y,Z]$ are complete in blue, $R^1:=[W,Y]$ and
$R^2:=[X,Z]$ are complete in red, and $G^1:=[W,Z]$ and $G^2:=[X,Y]$ are complete
in green.

\item there exists a partition $\{W, X, Y, Z\}$ of $[n]$ with $X,Y,Z$ non-empty
such that $B:=W\cup X\cup Y$ is connected in blue, $R:=W\cup X\cup Z$ is
connected in red, and $G:=W\cup Y\cup Z$ is connected in green.  Furthermore,
$[X,Y]$ is complete in blue, $[X,Z]$ is complete in red, and $[Y,Z]$ is complete
in green, whereas no edge in $[W,X]$ is green, no edge in $[W,Y]$ is red, and no
edge in $[W,Z]$ is blue.
\end{enumerate}
\end{lemma}

\begin{proof}
Suppose $B$ is a maximal monochromatic, say blue, connected subgraph and set $U=V(K)\setminus B$.
If $U=\emptyset$ then we are in case (i); so suppose not.  Note that all edges from $B$ to $U$ are either red or green.  Let $R$ be a maximal, say red, component which intersects both $B$ and $U$. By the maximality of $B$, we have $B\setminus R\neq \emptyset$. 

First suppose $U\setminus R\neq \emptyset$.  In this case, both $[B\cap R,
U\setminus R]$ and $[B\setminus R, U\cap R]$ are complete in green.  This
implies $[B\cap R, U\cap R]$ and $[B\setminus R, U\setminus R]$ are complete in
red and $[B\cap R, B\setminus R]$ and $[U\cap R, U\setminus R]$ are complete in
blue.  So we are in case (ii), setting $W:=B\cap R$, $X:=B\setminus R$, $Y:=U\cap R$,
and $Z:=U\setminus R$.

Finally, suppose $U\setminus R=\emptyset$.  In this case $[B\setminus R, U]$ is
complete in green, so there is a maximal green component $G$ containing $U\cup
(B\setminus R)$.  Then we are in case (iii), setting $W:=B\cap R\cap G$, $X:=B\setminus
G$, $Y:=B\setminus R$, and $Z:=U$. 
\end{proof}

\begin{figure}[ht]
\centering
\includegraphics[scale=0.85]{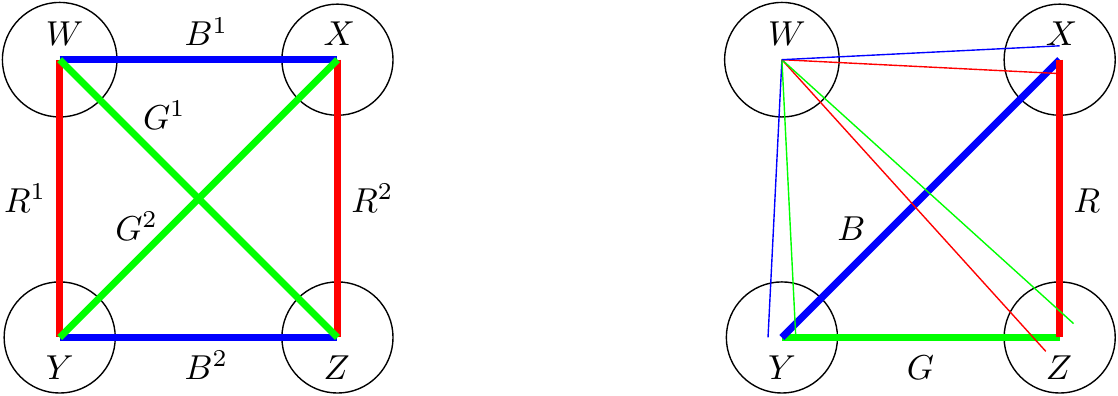}
\caption{Colorings of type (ii) and (iii) respectively.}\label{fig:types}
\end{figure}

First we give the short proof of the 2-color case.

\begin{proof}[Proof of Theorem \ref{infTree2}]
For each $n\in \NN$, the graph induced by $[n]$ is connected in either red or blue.  By pigeonhole, there exists, without loss of generality, an increasing sequence $n_1<n_2<\dots$ such that the graph induced by each $[n_i]$ is connected in red.  So there exists a red connected subgraph $R$ such that $d_s(R)=1$ (see the comment following Definition \ref{def:sud}).
\end{proof}

Now we prove the three color case.  Our original proof was quite complicated and we thank one of the referees for providing the following simpler proof.

\begin{proof}[Proof of Theorem \ref{infTree3}]
For each $n\in \NN$, say that $n$ is of type (i), (ii), or (iii) if Lemma \ref{lem:3col} (i), (ii), or (iii) holds respectively in the graph induced on $[n]$.

If there are infinitely many $n$ of type (i), then by pigeonhole, there is a
subsequence $n_1<n_2<\cdots$ such that the graph induced by each $[n_i]$ is
connected in the same color, and so we have a monochromatic connected subgraph
with strong upper density 1.  We may thus assume that there are only finitely
many $n$ of type (i).  Fix $n_0$ so that every $n\ge n_0$ is of type (ii) or
(iii).

Suppose that there is some $n\ge n_0$ of type (ii).  Let $\{W,X,Y,Z\}$ be the
partition of $[n]$ as described in Lemma \ref{lem:3col} case (ii), and define
$B^j$, $R^j$, and $G^j$ for $j = 1,2$ as in that case.  Consider the vertex $v =
n+1$.  This vertex cannot have blue edges to both $B^1$ and $B^2$, or red edges
to both $R^1$ and $R^2$, or green edges to both $G^1$ and $G^2$; otherwise $n+1$
would be of type (i).  Now note that $v$ must have edges of the same color to
two of the sets $W,X,Y,Z$; without loss of generality, $v$ has blue edges to
both $W$ and $X$.  Then every edge from $v$ to $B^2 = Y\cup Z$ is red or green;
without loss of generality, $v$ has a red edge to $Y$.  It follows that every
edge between $v$ and $Z$ is green, and this in turn implies that every edge
between $v$ and $Y$ is red, and every edge betwen $v$ and $X$ is blue.  Hence
$\{W\cup\{v\},X,Y,Z\}$ is a partition of $[n+1]$ meeting the requirements of
Lemma \ref{lem:3col} case (ii).  Applying this argument repeatedly, it
follows that every $m\ge n$ is of type (ii), and moreover we may asume that the
partitions $\{W_m,X_m,Y_m,Z_m\}$ of $[m]$ given by Lemma \ref{lem:3col} satisfy
$W_m\subseteq W_{m+1}$, $X_m\subseteq X_{m+1}$, $Y_m \subseteq Y_{m+1}$, and
$Z_m \subseteq Z_{m+1}$.  Thus, every color has a monochromatic component with
strong upper density at least $1/2$.

Now we may assume that every $n\ge n_0$ is of type (iii).  For each $n\ge n_0$,
fix a partition $\{W_n,X_n,Y_n,Z_n\}$ of $[n]$ as in Lemma \ref{lem:3col} case
(iii).  We define $R_n$, $B_n$, and $G_n$ as in that case (see Figure \ref{fig:types}), and note that at least two of the sets $R_n$, $B_n$ and $G_n$ are of order at least $n/2$.  By
pigeonhole we may find a sequence $n_1 < n_2 < \cdots$ such that, without loss
of generality, $|R_{n_i}|,|B_{n_i}|\ge {n_i}/2$.  By Ramsey's theorem, we may
pass to a subsequence where one of the following happens for all $i < j$;
$R_{n_i}\cap R_{n_j}\neq \emptyset$; $B_{n_i}\cap B_{n_j}\neq \emptyset$; there
is a red edge between $R_{n_i}$ and $R_{n_j}$; there is a blue edge between
$B_{n_i}$ and $B_{n_j}$; or $R_{n_i}\cap R_{n_j} = B_{n_i}\cap B_{n_j} =
\emptyset$ and every edge in $[R_{n_i}\cap B_{n_i},R_{n_j}\cap B_{n_j}]$ is
green.  All but the last of these cases imply that we have a red or blue
connected subgraph with strong upper density at least $1/2$, so we may assume
the last case holds.  Note now that every edge between $X_{n_i}$ and $X_{n_j}$
is green, so if $|X_{n_i}|\ge n_i/2$ for infinitely-many $i$ we are done; so  
we may assume that $|X_{n_i}| < n_i/2$, and hence $|G_{n_i}| \ge n_i/2$, for all
$i$.

If there are infinitely-many $i$ for which $W_{n_i}$ is nonempty, then
by another application of Ramsey's theorem, we may pass to a subsequence where
for all $i < j$, $W_{n_i}\cap W_{n_j} \neq\emptyset$, or some color, say red,
there is a red edge between $W_{n_i}$ and $W_{n_j}$; in either case we
get a monochromatic connected component with strong upper density at least
$1/2$.

So suppose that $W_{n_i}$ is empty for all large enough $i$.  Recall that all
edges in $[X_{n_i},Y_{n_i}]$ are blue, all edges in $[X_{n_i},Z_{n_i}]$ are red,
and all edges in $[Y_{n_i},Z_{n_i}]$ are green.  Furthermore, since all of the
components $R_{n_i}$, $B_{n_i}$, and $G_{n_i}$ are of order at least $n_i/2$, we
may assume by another application of Ramsey's theorem that for all $i < j$, all edges in $[X_{n_i}, X_{n_j}]$ are green, all edges in $[Y_{n_i},Y_{n_j}]$ are
red, and all edges in $[Z_{n_i},Z_{n_j}]$ are blue.  It follows that $X_{n_i}
\subseteq X_{n_{i+1}}$, $Y_{n_i}\subseteq Y_{n_{i+1}}$, and $Z_{n_i}\subseteq
Z_{n_{i+1}}$.  Now we see that each color has a component with strong upper
density at least $1/2$.   
\end{proof}

\section{0-dense graphs}\label{sec:0dense}

We are interested in which graphs $G$ have the property that in every $r$-coloring of the edges of $K_\NN$, there exists a monochromatic copy of $G$ with positive upper density.  We give some partial answers to this question in the undirected setting as well as the directed setting.  

\begin{proposition}\label{prop:relden}
Let $h : \NN\to\NN$, where $h(n)\to\infty$.  There exists a $2$-coloring of $K_\NN$ such that if $G$ is a graph in which the largest independent set has size $i<\infty$, then the vertices of any monochromatic copy of $G$ must eventually grow faster than $h$.
In particular, any graph with bounded independence number is 0-dense.
\end{proposition}

\begin{proof}
Let $\NN = A_1\cup A_2 \cup \cdots \cup A_n \cup \cdots$ be a partition of $\NN$
into consecutive intervals, such that $|A_n| = nh(n)$.  Color the edges inside
the intervals blue, and the edges between intervals red.  Let $G$ be an infinite graph
with independence number $i < \infty$.  Then if $V$ is the vertex set of a
monochromatic (necessarily red) copy of $G$, it follows that $|V\cap A_n| \le
i$ for each $n$, and from this that the vertices of $V$ eventually grow faster
than $h$.
\end{proof}

Let $\CT^+$ be the transitive tournament on $\NN$ such that for all $i,j\in \NN$, $(i,j)\in \CT^+$ if and only if $i<j$.
Let $\CT^-$ be the transitive tournament on $\NN$ such that for all $i,j\in \NN$, $(i,j)\in \CT^-$ if and only if $i>j$.  Given a directed graph $D$, $v\in V(D)$, and $S\subseteq V(D)$, define $N^+(v, S)=\{u\in S: (v,u)\in E(D)\}$ and $N^-(v, S)=\{u\in S: (u,v)\in E(D)\}$.  Then define $d^+(v,S)=|N^+(v, S)|$ and $d^-(v,S)=|N^-(v, S)|$.

  \begin{definition}
    We define two classes of countably infinite directed graphs, $\SF^+$ and $\SF^-$, as follows.  $\SF^\pm$ (we use the
    symbol $\pm$ to define them simultaneously) is defined to be the class of countably infinite graphs $G$ such that
    there is a partition $V(G) = A\cup B$ into infinite sets satisfying
    \begin{enumerate}
      \item  $A\cong \CT^\pm$ and $B$ is independent,
      \item  for each $a\in A$, $d^\mp(a,B)$ is finite,
      \item  for each $b\in B$, $d^\mp(b,A) = 0$.
    \end{enumerate} 
  \end{definition}
  In the following, a \emph{filter} is a collection $\SC$ of subsets of
  $\NN$ which is closed under finite intersections and supersets, and does not
  contain the empty set.  A filter $\SU$ is an \emph{ultrafilter} if for every
  $X\subseteq\NN$ either $X\in\SU$ or $\NN\sm X\in\SU$; equivalently, $\SU$ is
  an ultrafilter if $\SU$ is maximal among all filters with respect to
  inclusion.  An ultrafilter $\SU$ is \emph{nonprincipal} if every set in $\SU$ is
  infinite.  The existence of a nonprincipal ultrafilter follows easily from
  Zorn's lemma; indeed, Zorn's lemma implies that for any filter $\SC$ there is
  an ultrafilter $\SU$ containing $\SC$, so to guarantee that $\SU$ is
  nonprincipal we may simply take $\SC$ to be the filter of cofinite sets.  We note here that if $\SU$ is an
  ultrafilter, $A\in\SU$, and $A = A_1\cup\cdots\cup A_n$ is a partition of $A$,
  then exactly one $A_i$ is in $\SU$.

  \begin{theorem}\label{thm:ultratournament}
    Let $\SU$ be a nonprincipal ultrafilter.  Let $T$ be a tournament on $\NN$, and let a $k$-coloring of $T$
    be given.  Then the following is satisfied for at least one $s\in\{+,-\}$: for every $G\in \SF^s$ there is
    a $V\in\SU$ such that there is a monochromatic isomorphic copy of $G$ in $T$ with vertex set $V$.
  \end{theorem}
  
  \begin{proof}
    Fix $T$ and a $k$-coloring of $T$.  For each $v\in \NN$ we have either $N^+_T(v)\in\SU$ or $N^-_T(v)\in\SU$; put
    $v\in S^+$ in the first case and $v\in S^-$ in the second case.  Then the two sets $S^\pm$ partition $\NN$ and hence
    one must be in $\SU$.  We will show that, if $S^\pm\in\SU$, then for all $G\in \SF^\pm$ there is a monochromatic
    copy of $G$ in $\SU$.  We will do this in the case $S^+\in\SU$; the other case is symmetric.

    So we assume $S^+\in\SU$.  As above we partition $S^+$ into sets $S_i$, $i = 1,\ldots,k$, so that $v\in S_i$ if and
    only if $N_i^+(v)\in\SU$.  Since these sets partition $S^+$, one must be in $\SU$; say without loss of generality
    that $S_1\in\SU$.

    Fix $G\in \SF^+$, and let $V(G) = A\cup B$ be a partition as in the definition of $\SF^+$.  Let $a_i$ and $b_i$
    ($i\in\NN$) enumerate the elements of $A$ and $B$ respectively, and choose a function $f : \NN\to\NN$ satisfying
    \[
      N^-(a_i) \cap B \subseteq \set{b_j}{j \le f(i)}
    \]
    (To avoid degenerate cases we also require $f(i)\to\infty$ as $i\to\infty$.)

    Now we recursively choose vertices $x_1,x_2,\ldots$ and $y_1,y_2,\ldots$ from $S_1$ so that
    \begin{enumerate}
      \item  for all $i < j$, $(x_i,x_j)$ has color $1$,
      \item  for all $i$ and $j\le f(i)$, $(y_j,x_i)$ has color $1$, and
      \item  the sets $X = \set{x_i}{i\in\NN}$ and $Y = \set{y_i}{i\in\NN}$ partition $S_1$.
    \end{enumerate}
    To begin the recursion, let $y_1,\ldots,y_{f(1)}$ be the first $f(1)$-many vertices of $S_1$.  Let $A_i =
    N_1^+(y_i)$ for each $i\le f(1)$; then each $A_i$ is in $\SU$, and hence there are infinitely-many $x$ in $A_1\cap
    \cdots \cap A_{f(1)}\cap S_1$; let $x_1$ be the minimal such $x$ greater than $y_1,\ldots,y_{f(1)}$.

    Now suppose we have chosen $x_1,\ldots,x_n$ and $y_1,\ldots,y_{f(n)}$.  Let $y_{f(n)+1},\ldots,y_{f(n+1)}$ be the next
    $(f(n+1) - f(n))$-many vertices of $S_1$ not so far listed.  As in the base case, if we let $A_i = N_1^+(y_i)$ for
    each $i\le f(n+1)$ and $B_j = N_1^+(x_j)$ for each $j\le n$, then the set
    \[
      S_1\cap A_1\cap \cdots \cap A_{f(n+1)}\cap B_1\cap \cdots \cap B_n
    \]
    is in $\SU$, and hence we may choose $x_{n+1}$ to be some element of this set larger than all $x_i$'s and $y_i$'s
    listed so far.

    Given sets $X$ and $Y$ as described above, we can clearly find a monochromatic copy of $G$ with vertex set
    $S_1$.  This completes the proof.
  \end{proof}

An ultrafilter $\SU$ is called \emph{Ramsey} if for every coloring of $K_\NN$ in
finitely-many colors, there is an $A\in\SU$ such that the complete graph on $A$
is monochromatic.  (See~\cite{KT} Chapter 17 for other properties of Ramsey
ultrafilters.)  Hence if $\SU$ is Ramsey, then any infinite graph satisfies the
conclusions of Theorem \ref{thm:ultratournament}, that is, if $G$ is any
infinite graph and $K_\NN$ has been $k$-colored then there is a $V\in\SU$ for
which there is a monochromatic copy of $G$ with vertex set $V$.  However, it is
known that the existence of Ramsey ultrafilters is independent of ZFC (see
\cite{W}).  Moreover, every Ramsey ultrafilter necessarily contains sets whose
elements grow arbitrarily fast; this is implied by the results of~\cite{E} but
can also be proven directly.

On the other hand, it is well known that the subsets of $\NN$ of density 1 form
a filter and that any filter can be extended to an ultrafilter (by Zorn's
lemma).  So there exists an ultrafilter $\SU$ such that every $A\in \SU$ has
positive upper density (see \cite{KT} Chapter 17, Problem 7).  Applying Theorem
\ref{thm:ultratournament} with this special ultrafilter $\SU$, we get the
following corollaries;

\begin{corollary}\label{cor:infswitch}
Let $P = (p_1,p_2,\ldots)$ be a directed path with infinitely many switches.
Then for any $r$-coloring of a tournament on $\NN$ there is a monochromatic copy
of $P$ with positive upper density.
\end{corollary}

Say that a graph $G$ is \emph{half-locally finite} if $G$ contains an infinite
independent set $B$ and that for all $v\in V(G)\setminus B$, $d(v, B)$ is
finite.  (Note that $G[V(G)\setminus B]$ could induce a complete graph for
instance).  Equivalently, $G$ is half-locally finite if $G\subseteq F\in \SF$,
where $\SF$ is the family of undirected graphs obtained by removing the
orientations from say $\SF^+$.  

\begin{corollary}\label{cor:AB}
Let $G$ be a half-locally finite graph.  Then for any $r$-coloring of $K_\NN$, there is a monochromatic copy of $G$ with positive upper density.
\end{corollary}


\section{Open problems and concluding remarks}\label{sec:conclusion}

\subsection{Paths}

We conjecture the following extension of Theorem \ref{mainThm1}.

\begin{conjecture}
Every 2-coloring of $K_\NN$ contains a monochromatic path of upper density at least $8/9$.  If true, this is best possible by Example \ref{ex:2colorupper}.
\end{conjecture}

For finite graphs, Gy\'arf\'as, Ruszink\'o, S\'ark\"ozy, and Szemer\'edi \cite{GyRSS3} proved that for sufficiently large $n$, every $3$-coloring of $K_n$ contains a monochromatic path of length at least $n/2$.  We conjecture the following analog for infinite graphs. 

\begin{conjecture}
Every $3$-coloring of $K_\NN$ contains a monochromatic path of upper density at
least $1/2$. If true, this is best possible by Corollary \ref{ex:affine}.
\end{conjecture}

\subsection{Two-way infinite paths, i.e.\ cycles}

While we do not study this variant, it is likely that our result for strong upper density of paths extends to a two-way infinite version, which would be an analog of cycles for infinite graphs.  However, since a result for two-way infinite paths does not imply a result for one-way infinite paths, we did not explore this possibility.  On the other hand, for ordinary upper density, it is clear that our construction can be modified to give a two-way infinite path with upper density at least $3/4$.  

\subsection{Dense paths in both colors}

We mention a variant of the problem which seems not to have been studied even for finite graphs.  That is, it should be possible to extend Theorem \ref{thm:gg} to say that in every 2-coloring of $K_n$, there exists a red path $P_R$ with more than $\alpha n\geq 2n/3$ vertices and a blue path $P_B$ with more than $2(1-\alpha)n$ vertices (even stronger would be the additional conclusion that $V(P_R)\cup V(P_B)=V(K_n)$).  Again, it seems possible that such a result would extend to upper density and strong upper density versions, but we do not explore that possibility here.  

\subsection{Directed paths}

For consistently oriented paths (equivalently paths with finitely many
switches), the only case we were not able to settle is the following.\footnote{While this paper was under review, a 2-coloring of $\vec{K}_\NN$ in which every consistently oriented monochromatic path has upper density 0 was independently provided by Corsten \cite{Cor} and Guggiari \cite{Gug}.}

\begin{problem}\label{prob:directedpath2}
Does there exist a 2-coloring of $\vec{K}_\NN$ in which every consistently oriented monochromatic path has upper density 0?
\end{problem}

We showed that if a path has infinitely many switches it is positively-dense, but is there an absolute lower bound?  If we quantify how the switches are distributed, can we say something more precise?

\begin{problem}
Is it true that for every $\ep>0$, there exists a path $P$ with infinitely many switches and a 2-coloring of $\vec{K}_n$ such that every monochromatic copy of $P$ has upper density less than $\ep$?  

Let $f$ be an increasing function from $\NN$ to $\NN$ and let
$P_f=v_1v_2v_3\dots$ be a directed path such that $v_{f(i)}$ is a switch.
Determine a lower bound on the upper density of monochromatic copy of $P$ in
$\vec{K}_\NN$.  For instance if $f(i)=i$ for all $i\in \NN$, then $P_f$ is an
anti-directed path and a result of Rado from \cite{R} implies that every
2-coloring of $\vec{K}_\NN$ contains a monochromatic anti-directed path of upper
density at least $1/2$.   
\end{problem}

\subsection{Connected subgraphs}

We conjecture the following infinite analog of Theorem \ref{thm:Gy}, for which we provided the proof when $r\in \{2,3\}$.  

\begin{conjecture}
For $r\geq 4$, every $r$-coloring of $K_\NN$ contains a monochromatic connected subgraph $T$ with $\bar{d}(T)\geq \bar{d}_S(T)\geq 1/(r-1)$.  If true, this is best possible when $r-1$ is a prime power by Example \ref{ex:affine}.
\end{conjecture}

\subsection{Relative density}

We showed that consistently oriented paths are $(0,3)$-dense and graphs with bounded independence number are $0$-dense, and in both cases we showed that no bound can be put on the growth of the vertices in a monochromatic copy.  Is there a $0$-dense graph for which some bound can be placed on the growth?

\begin{problem}
Does there exist a graph $G$ and a 2-coloring of $K_\NN$ such that for every monochromatic copy of $G$, we have $\bar{d}(G)=0$, but there exists an increasing unbounded function $h$ such that there exists a monochromatic copy of $G$ with $\bar{d}_h(G)=\limsup_{n\to \infty}\frac{h(|V(G)\cap [n]|)}{h(n)}>0$?
\end{problem}

\subsection{0-dense graphs}

Chv\'atal, R\"odl, Szemer\'edi, and Trotter \cite{CRST} proved that for all $\Delta\in \NN$, there exists $0<c\leq 1$ such that if $G$ is a graph on at most $cn$ vertices with maximum degree at most $\Delta$, then in every $2$-coloring of $K_n$, there exists a monochromatic copy of $G$.  We conjecture the following infinite analog of their result.  

\begin{conjecture}
For all $\Delta\in \NN$ there exists $c>0$ such that for every infinite graph $G$ with maximum degree at most $\Delta$ and every 2-coloring of $K_\NN$, there exists a monochromatic copy of $G$ with $\bar{d}(G)\geq c$.
\end{conjecture}

We showed in Corollary \ref{cor:AB} that if a graph $G$ is half-locally finite,
then it is positively-dense, but is there an absolute lower bound?  If so, this
would strongly imply the previous conjecture.  

\begin{problem}
For all $\ep>0$, does there exists a half-locally finite graph $G$ and a 2-coloring of $K_\NN$ such that every monochromatic copy of $G$ has upper density less than $\ep$? If so, what about a locally finite graph $G$?
\end{problem}

We show that graphs with bounded independence number are 0-dense and half-locally finite graphs are positively-dense.  It would be nice to have a complete characterization.

\begin{problem}
Completely characterize which graphs $G$ are positively-dense. 
\end{problem}

\subsection{Lower density}

\begin{problem}
Let $G$ be a graph.  Is it true that if there exists a $2$-coloring of $K_\NN$ such that for every monochromatic copy of $G$ we have $\bar{d}(G)<1$, then there exists $2$-coloring of $K_\NN$ such that for every monochromatic copy of $G$ we have $\ubar{d}(G)=0$?
\end{problem}

\subsection{Density results from partitioning results}

In a recent paper, Elekes, D. Soukup, L. Soukup, and Szentmikl\'ossy \cite{ESSS}, proved a variety of results about partitioning infinite (hyper)graphs.  A number of these results have consequences for the density Ramsey problem studied in this paper.  They prove the following three results:

(i) In any $r$-coloring of $K_\NN^{k}$ (the complete $k$-uniform hypergraph on $\NN$), there exists a partition of $\NN$ into at most $r$ monochromatic tight paths, and consequently there exists a monochromatic tight path with upper density at least $1/r$.  The analogous problem for finite hypergraphs was studied in \cite{GySSz}, \cite{HLPRRSS} for loose cycles, and \cite{HLPRRS} for tight cycles.

(ii) In any $2$-coloring of $K_\NN$, there exists a partition of $\NN$ into at most four monochromatic squares of paths and thus there exists a monochromatic square path of upper density at least $1/4$.  
More generally, they showed that in any $r$-coloring of $K_\NN$, there exists a partition of $\NN$ into at most $r^{(k-1)r+1}$ monochromatic $k$th powers of paths and a finite set, and thus there exists a monochromatic $k$th power of a path with upper density at least $\frac{1}{r^{(k-1)r+1}}$.  The analogous problem for squares of paths in 2-colored finite complete graphs was considered in \cite{ARS}.

(iii) In any $r$-coloring of the infinite random graph on $\NN$ there exists a partition of $\NN$ into at most $r$ monochromatic paths and thus some monochromatic path has upper density at least $1/r$. 
The analogous problem for finite random graphs was considered in \cite{Letz}.

It would be interesting to improve any of these density results (as Erd\H{o}s and Galvin improved the corollary of Rado's result) and to consider the corresponding strong upper density versions.

\section{Acknowledgements}

We thank the referees for their detailed reports, especially for providing a simplification of our original proof of Theorem \ref{infTree3}.

\section*{Appendix A: Erd\H{o}s-Galvin examples}

The following two examples appear as Corollary 1.6 and Corollary 1.9 in \cite{EG}.  As they are stated in more general terms in their paper, we reproduce them here in less general terms for the convenience of the reader.

\begin{example}\label{ex:2colorstrongupper}
  There exists a $2$-coloring of $K_\NN$ such that for all monochromatic paths $P$, $\bar{d}_s(P)\leq 2/3$.
\end{example}

\begin{proof}
  Let $A$ be the set of $n\in\NN$ divisible by $3$, and $B$ its complement.  Color every edge from $A$ to $B$ red; and color every edge from $A$ to $A$, and from $B$ to $B$, blue.  There is clearly a
  blue path which covers all of $B$ in order, and this path has strong upper density, and upper density, equal to $2/3$.  
  Now suppose $P = (p_0,p_1,\ldots,p_n,\ldots)$ is any red path; say, without loss of generality, that $p_0\in A$.  Then for
  every $n$, $p_{2n}\in A$ and $p_{2n} \ge 3n$.  Hence the density of $p_0,\ldots,p_{2n}$ inside any interval that contains this initial segment must be at most $2/3$.
\end{proof}

\begin{example}\label{ex:2colorupper}
  There exists a $2$-coloring of $K_\NN$ such that for all monochromatic paths $P$, $\bar{d}(P)\leq 8/9$. 
\end{example}

\begin{proof}
  Let $A_0, A_1, A_2, \dots$ be a sequence of consecutive intervals in $\NN$, with $|A_n|=2^n$ for all $n\geq 0$.  If $n$ is even and $m\geq n$, color all edges between $A_n$ and $A_m$ red.  If $n$ is odd and $m\geq n$, color all edges between $A_n$ and $A_m$ blue (See Figure \ref{fig:89}).
\begin{figure}[ht]
\begin{center}
\includegraphics[scale=.9]{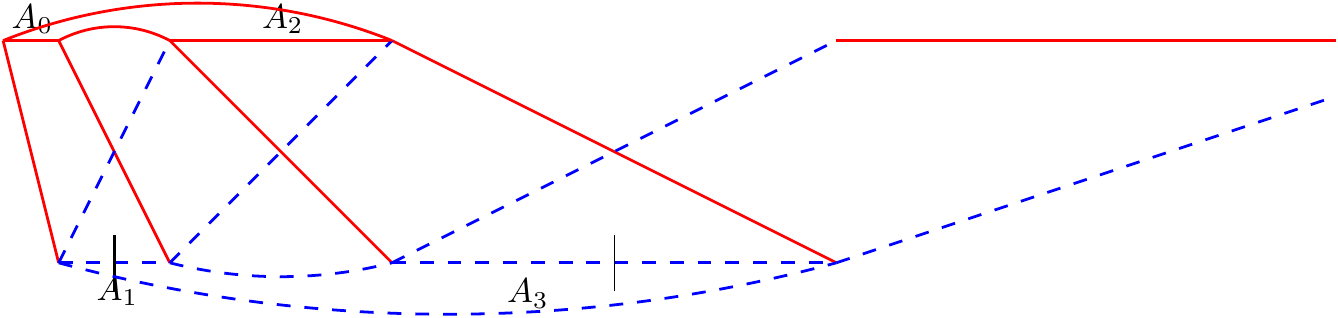}
\end{center}
\caption{An illustration of the example defined above.}\label{fig:89}
\end{figure}

By symmetry, we can consider the upper density of a red path $P$, which we can see is maximized by alternating between vertices in $A_0$ and $A_1$, using the earliest vertices in $A_1$ possible, until $A_0$ has filled up; then move  to $A_2$ and continue.  So the upper density of $P$ is
\begin{align*}
\lim_{n\to \infty} \frac{2\cdot 2^0 + 2\cdot 2^2 + \cdots + 2\cdot 2^{2n}}{2^0 + 2^1 + \cdots + 2^{2n} + 2^{2n}}=\lim_{n\to \infty} \frac{2\sum_{i=0}^n\frac{1}{4^i}}{1+\sum_{i=0}^n\frac{1}{2^i}}=\frac{2\cdot\frac{4}{3}}{1+2}=8/9.
\end{align*}

\end{proof}

\section*{Appendix B: Szemer\'edi's regularity lemma}\label{app:reg}

We will use the following version of Szemer\'edi's regularity lemma.  For a more detailed introduction, we refer the reader to \cite{KO} or \cite{KS}.

\begin{lemma}[2-colored regularity lemma -- degree form]\label{2colordegreeform}
For all $0< \ep\leq 1$ and positive integers $m$, there exists an $M = M(\ep,
m)$ such that for all 2-colored graphs $G$ on $n\geq M$ vertices and all
$d\in [0,1]$, there exists an integer $k$, a partition $\{V_0,
V_1,\dots, V_k\}$ of the vertex set $V(G)$, and a subgraph $G' \subseteq G$ with the following properties:
\begin{enumerate}
\item $|V_0| \le \ep n$
\item $m \le k \le M$ and $|V_1|=|V_2|=\dots=|V_k|$,
\item $d_{G'}(v)>d_G(v)-(2d+\ep)n$ for all $v\in V$,
\item $e(G'[V_i]) = 0$ for all $i \in [k]$,
\item for all $1 \le i < j \le k$, the pair $(V_i,V_j)$ is $\ep$-regular in $G_1'$ with a density either 0 or greater than $d$ and $\ep$-regular in $G_2'$ with a density either 0 or greater than $d$, where $E(G') = E(G_1') \cup E(G_2')$ is the induced 2-coloring of $G'$.
\end{enumerate}
\end{lemma}


\begin{definition}[$(\ep,d)$-reduced graph]\label{def:reduced}
Given a 2-colored graph $G$ and a partition $\{V_0, V_1, \dots, V_k\}$ satisfying the conclusion of Lemma \ref{2colordegreeform}, we define the $(\ep, d)$-reduced graph of $G$ to be the graph $\Gamma$ on vertex set $\{V_1, \dots, V_k\}$ such that $V_iV_j$ is an edge of $\Gamma$ if $G'[V_i, V_j]$ has density greater than $2d$.  For each $V_iV_j\in E(\Gamma)$, we assign color $1$ if $G'_1[V_i, V_j]$ has density greater than $d$ and color $2$ if $G'_2[V_i, V_j]$ has density greater than $d$ (note that since the total density is greater than $2d$ every edge must receive a color, but it need not be unique). 
\end{definition}

The following two Propositions are commonly used together with the regularity lemma (see  \cite[Proposition 42 and 43]{KO}).

\begin{proposition}\label{reduceddegree}
Let $0<2\ep\leq d\leq c/2$ and let $G$ be a graph on $n$ vertices with $\delta(G)\geq cn$. If $\Gamma$ is a $(\ep, d)$-reduced graph of $G$ obtained by applying Lemma \ref{2colordegreeform}, then $\delta(\Gamma)\geq (c-3d)k$.
\end{proposition}

\begin{proposition}\label{superreg}
Suppose that $4\ep<d\leq 1$ and that $C$ is a cycle in the reduced graph $\Gamma$ (including the case that $C$ is a single edge). Then every set $V_i$ contains a subset $V_i'\subseteq V_i$ of size $(1-2\ep)m$ such that $(V_i',V_j')_{G'}$ is $(2\ep, d-3\ep)$-super-regular for every edge $V_iV_j \in E(C)$.
\end{proposition}


\begin{customlem}{\ref{supercycles}}
For all $0<\eta\leq 1/4$ there exists $0<\frac{1}{n_0}\ll \frac{1}{m}\ll \ep< 4d\ll\sigma$ such that every $2$-colored graph $K_n$ on $n\geq n_0$ vertices has a monochromatic $(\ep, d, m)$-super-regular cycle $C_1$ on at least $(2/3-\eta)n$ vertices.
\end{customlem}

\begin{proof}
Apply Lemma \ref{2colordegreeform} to the 2-colored $K_n$ (with $\ep/2$ and $2d$) and then Proposition \ref{reduceddegree} to get 2-colored $(\ep/2, 2d)$-reduced graph $\Gamma$ on $k$ vertices with $\delta(\Gamma)\geq (1-6d)k\geq (3/4+\ep/2)k$.  Apply Theorem \ref{conmatch} to $\Gamma$ to get a monochromatic cycle on at least $(2/3-\ep/2)k$ vertices in $\Gamma$.  Now use Proposition \ref{superreg} to turn the cycle in $\Gamma$ into a monochromatic $(\ep, d, m)$-super-regular cycle in $G$ on at least $(2/3-\eta)n$ vertices.
\end{proof}

\begin{customlem}{\ref{blowupcor}}
Let $k\geq 2$ and $0<\frac{1}{m}\ll \ep, d$.  An $(\ep, d, m)$-super-regular cycle $\mathcal{C}$ of order $k$ has the property that for all $u,v\in V(C)$, there exists a $u,v$-path of length at least $(m-1)k$.
\end{customlem}

\begin{proof}
Let $V_1, V_2, \dots, V_k$ be a cyclic ordering (meaning $V_{1-1}=V_k$ and $V_{k+1}=V_1$) of the vertex set of $\mathcal{C}$ and let $u\in V_i$ and $v\in V_j$.  Let $U_{i+1}=N(u)\cap V_{i+1}$ and $U_{j-1}=N(v)\cap V_{j-1}$ and note that $|U_{i+1}|, |U_{j-1}|\geq dm$.  If each $V_iV_{i+1}$ were complete, we could embed a path starting in $U_{i+1}$ and ending in $U_{j-1}$ which misses at most one vertex from each set.  So we may use the blow-up lemma to embed a path $P\subseteq \mathcal{C}$ starting in $U_{i+1}$ and ending in $U_{j-1}$ which misses at most one vertex from each set.
\end{proof}

\end{document}